\documentclass[11pt,a4paper]{amsart}
\usepackage[utf8]{inputenc}
\usepackage{thmtools}
\usepackage{thm-restate}

\usepackage{hyperref}

\usepackage{cleveref}
\usepackage{amsmath,amsthm}
\usepackage{fullpage,stmaryrd}
\usepackage{amsfonts}
\usepackage{amssymb}
\usepackage{mathtools} 
\DeclarePairedDelimiter{\abs}{\lvert}{\rvert}
\DeclarePairedDelimiter{\nor}{\lVert}{\rVert}

\newcommand\dens[1]{\mathrm{d}(#1)}
\newcommand\cD{\mathcal{D}}

\newcommand\cU{\mathcal{U}}
\newcommand\N{\mathbb{N}}
\newcommand\Q{\mathbb{Q}}
\newcommand\R{\mathbb{R}}
\newcommand\Z{\mathbb{Z}}
\newcommand\T{\mathbb{T}}
\newcommand\mt{>}
\newcommand\gt{>}

\newcommand\lt{<}

\newcommand\E{\mathbb{E}}

\newcommand\ldens[1]{\underline{\mathrm{d}}(#1)}
\newcommand\udens[1]{\overline{\mathrm{d}}(#1)}
\newtheorem{theorem}{Theorem}[section]
\newtheorem*{theorem*}{Theorem}

\newtheorem{lemma}[theorem]{Lemma}
\newtheorem{cor}[theorem]{Corollary}
\theoremstyle{definition}
\newtheorem{defn}{Definition}
\title{Realisability of simultaneous density constraints for sets of integers}
\author{Pierre-Yves Bienvenu}
\address{Institute of Discrete Mathematics and Geometry, TU Wien,
Wiedner Hauptstr. 8--10, A-1040 Wien, Austria} 
\email{pierre.bienvenu@tuwien.ac.at}
\begin{document}
\begin{abstract}
In this note,
we study the set $\cD$ of values of the quadruplet 
$(\ldens{A},\udens{A},\ldens{2A},\udens{2A})$
where $A\subset\N$ and $\ldens{\cdot},\udens{\cdot}$ denote the lower and upper asymptotic density, respectively.
Completing existing results on the topic,
we determine each of its six projections on coordinate planes, that is, the sets of possible
values of the six subpairs of the quadruplet.
Further, we show that this set $\cD$ has non empty interior, in particular 
has positive measure.
To do so, we use among others probabilistic and diophantine methods. Some auxiliary results pertaining to these methods may be of general interest.
\end{abstract}

\maketitle
\section{Introduction}
\label{sec:intro}
Given a set $A\subset\N$, let $2A=A+A=\{a+b:(a,b)\in A^2\}$ denote its sum set,
and similarly let $A-A=\{a-b:(a,b)\in A^2\}$ denote its difference set.
Let $\ldens{A}=\liminf \frac{\abs{A\cap [N]}}{N}$
and $\udens{A}=\limsup \frac{\abs{A\cap [N]}}{N}$
denote its upper and lower asymptotic density.
Here $[N]=\{1,\ldots,N\}$ for an integer $N\in\N$.
Whenever $\ldens{A}=\udens{A}$, we say that $A$ admits (or has, or possesses) a natural density,
and we denote by $\dens{A}$ this common value called the natural density.

For a set $A\subset\N$, we call density profile the quadruplet $p_A=(\ldens{A},\udens{A},\ldens{2A},\udens{2A})$
This note contributes to the determination of the set
$\cD=\{p_A : A\subset\N\}\subset[0,1]^4$.
Recently, a series of papers
studied the realisability of various upper and lower density constraints
for sets and sum sets
\cite{grekosIndI,grekosIndII,BH,LEONETTI_TRINGALI_2024}.
As an aside, let us point out that the study of density constraints for sets, product sets, and subset sums has also recently gained ground \cite{hennecartAustralasian,bettin}.
Several constructions have been devised to produce sets of prescribed densities.
\begin{itemize}
    \item Archimedean construction: this consists in taking sets of the
    form $A=\bigcup_n I_n$, where the $I_n$'s are disjoint (and even far apart) intervals of integers of  increasing lengths. These are the ``lacunary'' sets
    of \cite{bilu}, those satisfying $\limsup a_{n+1}/a_n\gt 0$.
    \item Non-archimedean or periodic construction: there we take sets of the form
    $A=F+q\cdot\N$, where $q\geq 2$ is an integer, $F\subset \N$ is a finite set (which may be assumed without loss of generality to be included in $\{1,2,\ldots, q\}$.
    \item Diophantine construction or Bohr sets: this consists in considering
    sets of the form $\{n\in\N : \lambda n\in A\}$ for some irrational number $\lambda\in\R/\Z$ and a subset $A$
    of the circle $\R/\Z$.
    This can be generalised to multidimensional Bohr sets, where $\lambda \in (\R/\Z)^m$
    and $A\subset (\R/\Z)^m$ for some $m\in\N$.
    \item Probabilistic construction: here we generate random sets following a suitable probability distribution.
\end{itemize}
The  first three are examples of sets with small doubling, i.e.
the density (upper or lower) of $2A$ is not much larger than that of $A$. The last one yields typically sets of large doubling.
Of course, these methods may be combined, as we show in this paper. Here we will use all but the periodic type.

The first theorem we state is a collection of existing results on the values of subpairs of the quadruplet $p_A$. Equivalently, it describes the projection of $\cD$
on some of the six coordinate planes.
For any pair $\{i,j\}\subset \{1,2,3,4\}$,
let $$\cD_{i,j}=\{(a,b)\mid \exists (x_1,x_2,x_3,x_4)\in \cD,\, x_i=a,\,x_j=b\}$$
be the projection of $\cD$ on to the planes of coordinate set $\{i,j\}$.
The first one, obtained in \cite[Theorem 1]{grekosIndII}, confirms that apart from the obvious constraints $\udens{A}\geq \ldens{A}$ and $\udens{2A}\geq \ldens{2A}$, the pairs
$(\ldens{A},\udens{A})$ and $(\ldens{2A},\udens{2A})$ are entirely free.
The second one, obtained by the author and Hennecart \cite[Theorem 1.1, Theorem 1.5]{BH}, asserts that the constraint
arising from Kneser's theorem \cite{kneser} is the only one on the pair $(\ldens{A},\ldens{2A})$. 
\begin{theorem}[See \cite{grekosIndII,BH}]
    \label{th:subpairs}
    \begin{enumerate}
        \item $\cD_{1,2}=\cD_{3,4}=\{(x,y)\in [0,1]^2:y\geq x\}$.
        \item $\cD_{1,3}$ and $\cD_{1,4}$ are both equal to 
the set of all pairs $(\alpha, \beta)\in [0,1]^2$ such that
$\beta\geq \min(1,2\alpha)$
or 
$\beta\in\Q\cap (0,1)$ has nonpositive dyadic valuation and
$\beta/2<\alpha\leq \beta/2+1/(2g_0)$,
where $g_0=\min\{g\in\N : g\beta\in\N\}$. 
    \end{enumerate}
\end{theorem}
The pair $(\ldens{A},\udens{2A})$, i.e. $\cD_{1,4}$, is not explicitly covered by
\cite{BH}, which deals with $\cD_{1,3}$. But
it is easy to see that \cite{BH} actually also shows that
 the set
$\{(\ldens{A},\udens{2A}) : A\subset\N\}$ is exactly the same set;
indeed, the fact that the conditions outlined above
are necessary follows again from Kneser's theorem (and the fact that $\udens{2A}\geq \ldens{2A})$,
and they are sufficient since the construction of \cite{BH}
actually provides a set $A\subset\N$ such that
$A$ and $2A$ admit a natural density.

The next theorem is new and supplies the two missing projections.
Let $f$ be the function $[0,1]\rightarrow [0,1]$ defined by $f(\alpha)=\min(3\alpha/2,(1+\alpha)/2)$.
   Let $\mathcal{U}=\{\alpha=(\alpha_1,\alpha_2)\in [0,1]^2\mid\alpha_2\geq f(\alpha_1)\}$.
\begin{theorem}
\label{th:newProj}
    \begin{enumerate}
        \item $\cD_{2,4}=\{\frac{\alpha}{k}:\alpha\in\mathcal{U};k\in\N\}$
        \item $\cD_{2,3}=[0,1]^2$
    \end{enumerate}
\end{theorem}
Finally, we have some partial results on the whole quadruplet $p_A$ for $A\subset\N$.
We show an explicit set of quadruplets $(a,b,c,d)$ defined by linear equations,
that is, a polyhedron, is included in $\cD$. 
\begin{restatable}{theorem}{quadruples}
\label{th:quadruples}
    The set of all quadruplets $(a,b,c,d)\in (0,1]^4$ satisfying the following system of inequalities
    
    $$\left\{
    \begin{array}{ccc}
        a    & \leq  & b \\
        c & \leq & d \\
        \min(2b,1) &\leq  &d \\
       \min( a+b,1) & \leq & c\\
        d &  < & 2c \\
    \end{array}
    \right.
   $$
   
is included in $\cD$.
\end{restatable}
The polyhedron of non-empty interior is obtained by replacing $\min(2b,1)$ by $2b$
and $\min( a+b,1)$.

The paper is organised as follows.
The next two sections are devoted to establishing general results on the density of sets and sumsets, first of lacunary ones, then
of random ones, that will be useful in the sequel.
We then prove Theorem \ref{th:newProj}
in Section \ref{sec:proj}.
Before proving Theorem \ref{th:quadruples} in Section \ref{sec:quad},
we recall Weyl's equidistribution theorem and obtain quantitative bounds for the error term in some cases in Section \ref{sec:dioph}; that Section also contains some new results
on the sum of sets of integers defined by diophantine constraints, which combines well with 
the probabilistic results.
Section \ref{sec:3dim} provides interesting consequences of Theorem \ref{th:quadruples} to certain three-dimensional projections of $\cD$, partially answering questions of Hegyv\'{a}ri, Hennecart and Pach \cite{hennecartAustralasian}.
Finally, we collect some further remarks on the set $\cD$
in Section \ref{sec:remarks}, explaining why this set can presumably not be characterised so cleanly as its two-dimensional projections, describing some complicated constraints
that arise when considering the full quadruples $p_A$ (or already subtriples thereof) and not only its subpairs.
\subsection*{Notation}
We use the notation $\llbracket a,b\rrbracket =[a,b]\cap\Z$ for any $a,b\in\R$.
Similarly $\llbracket a,b) =[a,b)\cap\Z$ and $(a,b\rrbracket =(a,b]\cap\Z$.
Also $[x]=\llbracket 1,x\rrbracket$.
The notation $|\cdot|$ may denote the cardinality of a set as well as the modulus of a complex number.
A probability is denoted by $\Pr$ and the mathematical expectation by $\E$.

\subsection*{Acknowledgement}
The  author is supported by FWF Grant I-5554.
\section{Results on the densities of lacunary sets}
\label{sec:lac}
The results in this section supply details on the ``archimedean'' construction of sets described above.
The first lemma shows that the lower density of a union of disjoint intervals is obtained by considering the finite density until gaps, whereas the upper one is obtained at the end of the full intervals, which is not too surprising, and is particularly helpful when the full intervals or the gaps get very large.
\begin{lemma}
\label{lm:1b}
    Let $x_n,y_n$ be increasing sequence of integers (in particular, tending to infinity). 
    Suppose $x_{n+1}>y_n>x_n$ for every $n\in\N$. Let $I_n=\llbracket x_n,y_n\rrbracket$ and $A=\bigcup_nI_n$.
    Then $\ldens{A}=\liminf \frac{|A\cap [x_n]|  }{x_n}$ and
    $\udens{A}=\limsup \frac{|A\cap [y_n]|  }{y_n}$.
    
    Two particular cases:
    \begin{enumerate}
        \item If $x_n/y_n$ tends
    to 0 and $x_{n+1}/y_n$ tends to the real number $\beta\geq 1$,
   then the lower density of $A$ is $1/\beta$, and its upper density is $1$.
   \item If $x_n/y_n$ tends to $1-\alpha \in (0,1)$ and $x_{n+1}/y_n$ tends to infinity,
then the lower density of $A$ is 0, and its upper density is $\alpha$.
    \end{enumerate}

\end{lemma}
\begin{proof}
    Let $x\in\N$. If $x\in A$, let $n\in\N$ be such that 
    $x\in I_n$. Then $|A\cap [x_n]|/x_n\leq |A\cap [x]|/x\leq |A\cap [y_n]|/y_n$.
    Otherwise, $x\notin A$ so let $n\in\N$ be such that
    $y_n\lt x\leq x_{n+1}$.
    Then $|A\cap [x_{n+1}]|/x_{n+1}\leq |A\cap [x]|/x\leq |A\cap [y_n]|/y_n$.
    As a result, $\ldens{A}=\liminf \frac{|A\cap [x_n]|  }{x_n}$ and
    $\udens{A}=\limsup \frac{|A\cap [y_n]|  }{y_n}$. This yields the first part of the lemma.

    Now we prove the two particular cases.
For the first one, we observe that $(y_{n-1}-x_{n-1})/x_n\leq \frac{|A\cap [x_n]|  }{x_n}\leq y_{n-1}/x_n$ for every $n\in\N$ and both the right hand side
and the left-hand side of this inequality tend to $1/\beta$.
Further, $\frac{|A\cap [y_n]|  }{y_n}\geq (y_n-x_n)/x_n$ which tends to 1.

For the second one, it suffices to note that
$\frac{|A\cap [x_n]|  }{x_n}\leq y_{n-1}/x_n\rightarrow 0$    
    and $|A\cap [y_n]|=y_n-x_n+O(y_{n-1})=\alpha y_n+o(y_n)$.\end{proof}

The next one is a criterium for a lacunary set to have (upper)
density 0.
\begin{lemma}
\label{lm:dens0}
    Let $a_n,b_n$ be two sequences of integers satisfying
    $a_n\lt b_n\lt a_{n+1}$ for all $n$ and
    $I_n=\llbracket a_n,b_n\rrbracket$.
    Assume $a_{n-1}=o(a_n)$.
    Let $I=\bigcup_n I_n$.
    Let $X_n\subset I_n$ for each $n$ and $X=\bigcup X_n$.
    Assume $\abs{X_n \cap (\min X_n+[0,y)}\leq \epsilon y$
    for all $y\geq \epsilon a_n$ for all $\epsilon \mt 0$,
    all $n$ large enough and all $y\in\N$ (equivalently all $y\leq b_n-a_n$).
    Then $\dens{X}=0$.
\end{lemma}
\begin{proof}
    We bound $|X\cap [t]|/t$ for large enough $t\in\N$.
    If $b_n \lt t\lt a_{n+1}$ for some $n\in\N$,
    then $|X\cap [t]|/t\leq |X\cap [b_n]|/b_n$.
    Therefore we may assume $a_n\lt t \lt b_n$ for some $n$.
    Write $t=a_n+y$. Then
\begin{equation}
        \label{eq:barycentre3}
        \frac{|X\cap [t]|}{t}=\frac{a_{n}}{t}\frac{|X\cap [a_{n}]|}{a_{n}}
+\frac{y}{t}\frac{|X_n\cap (a_{n},a_{n}+y]|}{y}.
    \end{equation}
    Let $\epsilon\gt 0$.
    We have $|X\cap [a_{n}]|=\abs{X_n\cap [b_{n-1}]}+O(x_{n-1})\leq \epsilon b_{n-1}/3+O(x_{n-1})\leq \epsilon a_n/2$ if $n$ is large enough, by hypothesis.
    If $y\leq \epsilon a_n/2$, then
    equation \eqref{eq:barycentre3} implies that $\frac{|X\cap [t]|}{t}\leq \frac{|X\cap [a_{n}]|}{a_{n}}+\epsilon/2\leq\epsilon $
Now suppose $y\geq \epsilon x_n$.
    By hypothesis, we have precisely
    $\frac{|X_n\cap (a_{n},a_{n}+y]|}{y}\leq \epsilon y/2$ for $n$ large enough.
    This implies that
    $|X\cap [t]|/t\leq \epsilon$, for every $\epsilon\gt 0$
    and $t$ large enough, i.e. $\udens{X}=0$ as desired.
\end{proof}
Here is a fairly simple corollary concerning sumsets that will be directly applicable.
\begin{cor}
\label{cor:dens0}
    Let $x_n,y_n$ be two sequences of integers satisfying
    $x_n\lt y_n\lt x_{n+1}$ for all $n$ and
    $I_n=\llbracket x_n,y_n\rrbracket$.
    Assume $x_{n-1}=o(x_n)$.
    Let $I=\bigcup_n I_n$.
    Let $A_n\subset I_n$ for each $n$ and $A=\bigcup A_n$.
    Assume further that for every $\epsilon\in (0,1)$, for every $n$ large enough and every $z\geq \epsilon x_n$,
    we have
    $\abs{2I_n\cap (2x_n,2x_n+z]\setminus 2A_n}\leq \epsilon z$ as well as 
    $\abs{(I_{n-1}+I_n)\cap (x_{n-1}+x_n,x_{n-1}+x_n+z]\setminus (A_n+A_{n-1})}\leq \epsilon z$.
    Then $\udens{2I\setminus 2A}=0$.
\end{cor}
\begin{proof}
    First, we remark that
    $2I=\bigcup_n J_n$
    where $J_n=I_n+\bigcup_{k\leq n}I_k$.
    Similarly  $2A=\bigcup_n B_n$
    where $B_n=A_n+\bigcup_{k\leq n}A_k$.
    Let $C_n=2A_n\cup (A_n+A_{n-1})\subset B_n$
    and similarly $K_n=2I_n\cup (I_n+I_{n-1})\subset J_n$.
    Then $J_n\setminus K_n\subset [x_n,x_n+x_{n-1}]$.
    Since $x_{n-1}=o(x_n)$, Lemma \ref{lm:1b} implies that
    $2I\setminus\bigcup_n K_n\subset \bigcup J_n\setminus K_n\subset \bigcup_n [x_n,x_n+x_{n-1}]$ has density 0.
    Thus we are left with showing that
    $\bigcup_n K_n\setminus 2A$ has density 0.
Since $\bigcup_n K_n\setminus 2A\subset \bigcup_n K_n\setminus C_n\subset \bigcup_n 2I_n\setminus 2A_n\cup \bigcup_n (I_n+I_{n-1})\setminus (A_n+A_{n-1})$,
it suffices to show that
$X=\bigcup_n X_n$ and $Y=\bigcup_n Y_n$ have density 0, where
$X_n=2I_n\setminus 2A_n$
and $Y_n=(I_n+I_{n-1})\setminus (A_n+A_{n-1})$.
Both follow directly from Lemma \ref{lm:dens0} thanks to the hypothesis, with $a_n=2x_n,b_n=2y_n$ in the case of $X$
and $a_n=x_n+x_{n-1},b_n=y_n+y_{n-1}$ in the case of $Y$.
\end{proof}

\section{Probabilistic lemmatas}
\label{sec:proba}
Our first result may sound a bit complicated, but is actually very easy to prove and will turn out handy.
\begin{lemma}
\label{lm:basicProba}
    Let $I,J$ be finite subsets of $\N$.
    Let $K\subset I+J$ and $r\in\N$.
    Let $p,q\in [0,1]$ be real numbers.
    Suppose that for all $n\in K$, there exists (at least) $r$ representations
    $n=a_i(n)+b_i(n)$ with $(a_i(n),b_i(n))\in A\times B$ for every $i\in [r]$.
    Suppose $A$ and $B$ are random subsets of $I$ and $J$ respectively such that
    $\Pr(a_i(n)\in A)=r$ and $\Pr(b_j(n)\in B)=q$ and additionally
    the $2r$ events of the form $\{a_i(n)\in A\}$ and $\{b_j(n)\in B\}$ for $i\in [r],j\in[r]$
    are mutually independent.
    Then with probability at least $1-e^{-rpq/2}$, we have
    $|A+B|\geq |K|(1-e^{-rpq/2})$.
\end{lemma}
\begin{proof}
    Let $n\in K$. For each $i\in [r]$, the probability of the event $\{a_i(n)\notin A\text{ or }b_i(n)\notin B\}$ is $1-pq$ by hypothesis.
    Since these $r$ events are pairwise independent by hypothesis, the probability of their intersection is $(1-pq)^r$. Thus $\Pr(n\notin A+B)\leq (1-pq)^r$.
    Therefore
    \begin{align*}
        \E[|K\setminus (A+B)|]&=\sum_{x\in K}\Pr(x\notin A+B)\\
        &\leq 
        \abs{K}(1-pq)^{r}\\
        &\leq \abs{K}e^{-rpq}
    \end{align*}
      By Markov's inequality, it follows that
    $$\Pr(|K\setminus (A+B)|\geq \abs{K}e^{-rpq/2})\leq e^{-rpq/2},$$
    which immediately yields the desired conclusion.
\end{proof}

Let us define the type of random construction
we want to consider in this article.

\begin{defn}
Let $B$ be a set of integers and let $p\in [0,1]$.
A Bernoulli variable $\xi$ of parameter $p$ is an integer variable defined by $\Pr(\xi=1)=p=1-\Pr(\xi=0)$.
A Bernoulli  subset $A$ of parameter $p$ of $B$
is a random set $\{n\in B:\xi_n=1\}$, where $(\xi_n)_{n\in A}$
is a Bernoulli process, i.e. a sequence of independent Bernoulli variables of parameter $p$.
\end{defn}
In view of Lemma \ref{lm:basicProba}, the number of representations of an element in a sum set will play a role, so we introduce some relevant definitions.
\begin{defn}
    Let $I,J\subset\Z$. For any $x\in\Z$, we define its
    representation count with respect to $I,J$ to be $R'_{I,J}(x)=\abs{\{(i,j)\in I\times J:i+j=x\}}$.
    Its restricted representation count (with respect to $I,J$) is
    $R_{I,J}(x)=\abs{\{\{i,j\}\subset\Z: (i,j)\in I\times J,\, i\neq j,\, i+j=x\}}$. The functions $R'_{I,J}$ and $R_{I,J}$ are called representation (resp. restricted representation) function with respect to $I,J$.

    When the two sets $I,J$ are evident from the context, they may be omitted.

    Note that as soon as $I$ and $J$ are disjoint,
    these two counts are equal. In contrast,
    when $I=J$ and $x$ is even, we have $R(x)=(R'(x)-1)/2$.
\end{defn}
Here is a simple result about the representation function of intervals.
\begin{lemma}
\label{lm:repres}
    Let $I,J$ be intervals of lengths $2\leq m\leq n$.
    For any $r\leq m$, we have $R'_{I,J}(x)\gt r$ for all but at most $2r$ of the elements $x\in I+J$ (the exceptions being the smallest $r$ and the largest $r$ ones).

    Further, $R_{I,J}(x)\geq r$
    for all but at most $4r$ elements of $I+J$.
    
\end{lemma}
\begin{proof}
    We prove the first statement.
    Without loss of generality, we may translate $I$ in such a way that
    $I+J=\llbracket 1,M \rrbracket$, where $M=m+n-1$.
    Then $R'(x)=\min(x,m,m+n-x)$ for any $x\in I+J$.
In particular $R'(x)\gt r$ for any $x \in I+J$
unless $x\in [r]$ or $m+n-x\in [r]$, which concludes.
The second statement follows from the first one since
$R_{I,J}\geq (R'_{I,J}-1)/2$.
\end{proof}
We state a simple consequence of Lemma \ref{lm:basicProba}.
 \begin{cor}
 \label{lm:proba}
     Let $I,J$ be intervals of lengths $2\leq m\leq n$.
     Let $A,B$ be independent Bernoulli subsets of $I$ and $J$ of parameter $p,q$ respectively.
     Let $\epsilon \in (0,1)$ and suppose
     $pqm\geq 16\epsilon^{-2}\log(M\epsilon^{-1})$.
Then for any $\epsilon\gt 0$, with probability $1-\epsilon$, whenever $y\in (\epsilon m,m+n)$,
we have
$\abs{(A+B)\cap (\min(I+J)+\llbracket 0,y-1\rrbracket)}\geq (1-\epsilon) y$
and
$\abs{(A+A)\cap (\min(2I)+\llbracket 0,y-1\rrbracket)}\geq 
(1-\epsilon) y.$
 \end{cor}
\begin{proof}
We first deal with $A+B$, the case of $A+A$ being very similar.
Let $y\in (\epsilon m,m+n)$.
Let $r=\min(m,\epsilon y/4)$, thus $r\geq \epsilon^2m/4$.

Then we apply Lemma \ref{lm:basicProba} with $K=\{x\in I+J:x\lt \min(I+J)+y,\, R'(x)\geq r\}$.
By Lemma \ref{lm:repres}, we know that that $\abs{K}\geq y-2r\geq y(1-\epsilon/2)$
since $r\leq \epsilon y/4$;
also  $$e^{-rpq/2}\leq e^{-\epsilon^2pqm/8}\leq\epsilon/M\leq \epsilon/2.$$
Thus with probability $1-e^{-rpq/2}\geq 1-\epsilon/M$,
we have $\abs{(A+B)\cap [y]}\geq \abs{K}(1-e^{-rpq/2})\geq y(1-\epsilon)$
since $e^{-rpq/2}\leq \epsilon/M\leq \epsilon/2$.
This is true for each single $y\in (\epsilon m,m+n)$,
so the probability that
$\abs{(A+B)\cap [y]}\geq  y(1-\epsilon)$ for every 
$y\in (\epsilon m,m+n)$ is, by the union bound,
at least $1-\epsilon$.

    To deal with $A+A$, we use almost the same proof, with $I=J$, $p=q$ and $m=n$. 
    This
    time 
    let $r=\min(\lfloor (n-1)/2\rfloor,\epsilon y/8)$, thus $r\geq \epsilon^2m/8$.
    Then we apply Lemma \ref{lm:basicProba} with $K=\{x\in 2I:x\lt \min(2I)+y,\, R(x)\geq r\}$.
    By Lemma \ref{lm:repres}, we know that that $\abs{K}\geq y-4r\geq y(1-\epsilon/2)$
From there on the proof runs just as above.
\end{proof}

\section{Projection on coordinate planes}
\label{sec:proj}

\subsection{Coordinate plane $\{2,4\}$}
Since every non empty set $B\subset\N$ may be written in the form
$B=x+k\cdot A$ where $x\geq 0$ and $k\geq 1$ are integers and $0\in A\subset\N$ and $\gcd(A)=1$,
it suffices, when determining $\cD$, to consider 
sets $A\subset\N$ which have this property or the weaker property
$\gcd(A-A)=1$, equivalently $A$ is not included in an arithmetic progression of common difference greater or equal to 2.
Let
$$\cD'=\{(\ldens{A},\udens{A},\ldens{2A},\udens{2A}) : A\subset\N,\gcd(A-A)=1\}.$$
Again we denote the six projections on coordinate planes by $\cD'_{i,j}$ for $(i,j)\subset \{1,2,3,4\}$.
Then
$\cD=\{\frac{1}{k}\alpha :\alpha\in \cD',\, k\in\N_{\geq 1}\}$.

Thus to determine 
$\cD_{2,4}$,
it suffices to determine 
$\cD'_{2,4}=\{(\udens{A},\udens{2A}) : A\subset\N, \gcd(A-A)=1\}$.
Recall from Section \ref{sec:intro} the following definitions.
   Let $f$ be the function $[0,1]\rightarrow [0,1]$ defined by $f(\alpha)=\min(3\alpha/2,(1+\alpha)/2)$.
   Let $\mathcal{U}=\{(\alpha_1,\alpha_2)\in [0,1]^2\mid\alpha_2\geq f(\alpha_1)\}$.
We will prove the following.
\begin{theorem}
\label{th:24}
We have $\cD'_{2,4}=\cU$.
\end{theorem}
It is clear that this implies Theorem \ref{th:newProj} (1).

Now by a theorem of Freiman,
for any $A\subset\N$ satisfying $\gcd(A-A)=1$, we have
$\udens{2A}\geq f(\udens{A})$.
Thus $\cD'_{2,4}\subset \cU$. This result of Freiman is a simple consequence of a finitary one, see for instance \cite{jin}[Lemma 1.1].

So there remains to prove $\cD'_{2,4}\supset \cU$.
Before that,
we present a well known construction yielding a  weaker result,
and upon tweaking the construction, we will obtain the full result.
\begin{lemma}
    \label{lm:afa}
    Whenever $\alpha\in [0,1]$, the pair
    $(\alpha,f(\alpha))$
    belongs to $\cD'_{2,4}$.
\end{lemma}
\begin{proof}    
When $\alpha=0$, the statement is obviously true so let
$\alpha\in (0,1]$.
Consider a sequence $T_n$ of integers tending so quickly to infinity that $T_n=o(T_{n+1})$.
For definiteness, we can select for instance
$T_n=2^{n^2}$.
Then we let
$B=\bigcup_n B_n$ where $B_n=[(1-\alpha)T_n,T_n]$.
Clearly $B$ is not contained in any arithmetic progression of common difference at least two
since it contains two consecutive elements.
Further $B$ has  upper asymptotic density $\alpha$ by Lemma \ref{lm:1b}, second particular case.
Also $B+B=\bigcup_n C_n$ where $C_n=B_n+\bigcup_{k\leq n}B_k$.
Note that
$$\llbracket (1-\alpha)T_{n-1},T_{n-1}\rrbracket = B_{n-1}\subset \bigcup_{k< n}B_k\subset \llbracket 0,T_{n-1}\rrbracket.$$
It follows that
\begin{multline}
\llbracket 2(1-\alpha)T_n,2T_n\rrbracket \cup \llbracket (1-\alpha)(T_n+T_{n-1}),T_n+T_{n-1}\rrbracket \subset
C_n\\
\subset [2(1-\alpha)T_n,2T_n]\cup [(1-\alpha)T_n,T_n+T_{n-1}].
\end{multline}
Let $D_n=\llbracket 2(1-\alpha)T_n,2T_n\rrbracket \cup \llbracket (1-\alpha)T_n,T_n\rrbracket$.
Let $E_n=C_n\Delta D_n$.
Then $E_n\subset \llbracket (1-\alpha)T_n,(1-\alpha)(T_n+T_{n-1})\rrbracket \cup \llbracket T-n,T_n+T_{n-1}\rrbracket$.
Let $E=\bigcup_n E_n$ and $D=\bigcup_n D_n$.
We have $|E\cap [T_{n}+T_{n-1}]|=O(T_{n-1})$ and $|E\cap [(1-\alpha)(T_{n}+T_{n-1})]|=O(T_{n-1})$; by Lemma \ref{lm:1a},
and since 
$T_{n-1}=o(T_n+T_{n-1})$
we infer $\udens{E}=0$.
It follows easily that
$\udens{B+B}=\udens{D}$.
Now assume $\alpha <1/2$, so that the intervals
$[2(1-\alpha)T_n,2T_n]$ and $[(1-\alpha)T_n,T_n+T_{n-1}]$ are disjoint.
We have
$|D\cap [T_n]|=\alpha T_n +O(T_{n-1})$ and
$|D\cap [2T_n]|=3\alpha T_n +O(T_{n-1})$. By Lemma \ref{lm:1a}, we have
$\udens{D}=3\alpha/2=f(\alpha)$
as desired.
Assume now $\alpha\geq 1/2$. Thus $D_n=\llbracket (1-\alpha)T_n,2T_n\rrbracket$.
By Lemma \ref{lm:1a} again, we have
$\udens{D}=(1+\alpha)/2=f(\alpha)$
as desired.
\end{proof}

To finish the proof of Theorem \ref{th:24}, we will use the probabilistic method.


\begin{proof}[Proof of Theorem \ref{th:24}]
Let $(\alpha_1,\alpha_2)\in [0,1]^2$ satisfy
$\alpha_2\geq f(\alpha_1)$.
We want to show $(\alpha_1,\alpha_2)\in \cD'_{2,4}$.
We may assume $\alpha_2\neq 0$ since we already know that
$(0,0)\in\cD'_{2,4}.$
Let $\alpha=f^{-1}(\alpha_2)$ (which is well defined and distinct from 0 since $f$ is a bijection from 
$[0,1]$ to itself) and consider $B=\bigcup_n B_n$ as defined in the proof of Lemma \ref{lm:afa}
with this $\alpha$, i.e.
$B_n=\llbracket (1-\alpha)T_n,T_n\rrbracket$ where $T_n=2^{n^2}$.
Further let $\beta=\alpha_1/\alpha$.
Note that $\beta\leq 1$ since $f$ is an increasing function.
Consider, for each $n$, a Bernoulli subset $A_n\subset B_n$
of parameter $\beta_n$, where $\beta_n=\beta$ if $\beta\neq 0$
and $\beta_n\rightarrow 0$ sufficiently slowly (at a rate to be determined later) otherwise.
Also assume that the random sets $A_n$ for $n\in\N$ are independent of each other.
Then let $A=\bigcup_n A_n$.
Then almost surely, $A$ has density $\alpha\beta=\alpha_1$ by the Law of large numbers, and $A$ lies in no arithmetic progression of common difference more than 1 because it almost surely contains two consecutive elements (by the Borel-Cantelli lemma for instance). If $\beta=0$, this requires
$\beta_n$ to decay sufficiently slowly; it suffices to have
$\sum_n\beta_n^2\rightarrow\infty$.
We assume $\beta_n=1/\sqrt{n}$ from now on.
Our goal is to show that $\udens{2A}=\udens{2B}$ almost surely.
Since $2A\subset 2B$, it suffices to show 
$\dens{2B\setminus 2A}=0$ almost surely.
For that, we will apply Corollary \ref{cor:dens0}.
Let $\epsilon_n=1/n^2$.
Note that $\beta_n^2\alpha T_{n-1}\geq 16\epsilon_n^{-2}\log(2T_n\epsilon_n^{-1})$ for every $n$ large enough.

Let $X_n=2B_n\setminus 2A_n$ and $Y_n=(B_{n-1}+B_n)\setminus (A_{n-1}+A_n)$.

Let $S_n=(1-\alpha)T_n$.
By Corollary \ref{lm:proba}, for every $n$,
with probability at least $1-\epsilon_n$,
we have $|X_n\cap (2S_n+[0,y-1])|\leq \epsilon_n y$
and $|Y_n\cap (S_{n-1}+S_n+[0,y-1])|\leq \epsilon_n y$ for every $y\geq \epsilon S_n$.
Since $\sum_n \epsilon_n$ converges, by Borel--Cantelli, with probability 1,
we have $|X_n\cap (2S_n+[0,y-1])|\leq \epsilon y$
and $|Y_n\cap (S_{n-1}+S_n+[0,y-1])|\leq \epsilon_n y$ for every $y\geq \epsilon T_n$ and every $n$ large enough simultaneously.
Thus we can apply Corollary \ref{cor:dens0}:
we infer that the upper asymptotic density of $A+A$ is almost surely just that of $B+B$, which is
$f(\alpha)=\alpha_2$ as was shown in the proof of Lemma \ref{lm:afa}.
In particular, there exists $A\subset\N$ such that $\gcd(A-A)=1$ and $(\udens{A},\udens{2A})=(\alpha_1,\alpha_2)$ as desired.
\end{proof}
\subsection{Coordinate plane 2,3}
Now there remains to determine $\cD_{2,3}$, i.e. the set of all possible values of
$(\udens{A},\ldens{2A})$. 
We will prove $\cD_{2,3}=[0,1]^2$.
First we show that $\cD_{2,3}$ contains $\{1\}\times [0,1]$.

\begin{lemma}
\label{lm:1a}
For any $\alpha\leq 1$, there exists a set $A\subset\N$
such that 
$(\udens{A},\ldens{2A})=(1,\alpha)$.
\end{lemma}
\begin{proof}
First consider the case where $\alpha \neq 0$.
    Consider the quickly increasing sequences of real numbers $(x_n)_{n\in\N},(y_n)_{n\in\N}$ defined
    by $x_n=2^{n^2}$ and $y_n=2^{(n+1)^2}/\beta $ for 
    each $n\in\N$, where 
    $\beta\geq 2$ is a real number to be determined later.
    Consider $A= \bigcup_n I_n$
    where $I_n=\llbracket x_n,y_n\rrbracket$.
    Here and anywhere else, we use the notation $\llbracket a,b\rrbracket =\N\cap [a,b]$.
    Lemma \ref{lm:1b} implies that $\udens{A}=1$ (and additionally
    $\ldens{A}=1/\beta$).
    Further, we can show that
    $\ldens{2A}=2/\beta$.
    Indeed,
    $2A=\bigcup_{n}K_n$
    where $K_n=I_n+\bigcup_{m\leq n}I_m$.
    It is clear that 
    $K_n\subset \llbracket x_n+1,2y_{n}\rrbracket$.
    In fact $K_n=\llbracket x_n+1,2y_n\rrbracket$ since $I_n+I_n\cup I_n+I_1=\llbracket x_n+1,2y_n\rrbracket$.
    By Lemma \ref{lm:1b}, this implies $\ldens{2A}= 2/\beta$.
    Thus taking $\beta=2/\alpha$, we conclude.

    Now assume $\alpha=0$.
    This time let
    $x_n=2^{(2n)^2}$ and $y_n=2^{(2n+1)^2}$.
    Consider $A= \bigcup_n I_n$
    where $I_n=\llbracket x_n,y_n\rrbracket$.
    Then arguing along similar lines as above,
    we get $\udens{A}=1$ and $\ldens{2A}=0$ as desired.
\end{proof}
To realise any $(x,y)\in [0,1]^2$, we will 
take a random subset from the construction above, in a very similar way as the derivation of Theorem \ref{th:24} from the construction of Lemma \ref{lm:afa}.
\begin{theorem}
\label{th:23}
    The set of all possible values of
$(\udens{A},\ldens{2A})$ as $A$ ranges among the subsets of $\N$ is $[0,1]^2$.
\end{theorem}
\begin{proof}
    Let $(a,b)\in [0,1]^2$.
    If $a=0$, let $(a_n)$ be a sequence of reals in $[0,1]$ slowly decaying to 0; if $a\neq 0$, let $a_n=a$ for every $a$.
   Then consider the set $A=\bigcup_n I_n$ constructed in Lemma \ref{lm:1a} for $\alpha=b$. 
   In particular $\ldens{2A}=b$.
   We use the same notation
   as in the proof of that lemma.
   We now construct a random set of the form
   $B=\bigcup_n B_n$
   where $B_n\subset I_n$ is a Bernoulli subset of $I_n$ of parameter $x_n$ for every $n$, and assume the random sets $B_n$ are pairwise independent. 
Then almost surely, this set has an upper density $a$ by the Lax of Large numbers.
Now we want to show $\ldens{2B}=\ldens{2A}=b$.
It suffices to show $\dens{2A\setminus 2B}=0$,
for which we will apply Corollary \ref{cor:dens0}.
To do so, we need to study
$X_n=2I_n\setminus 2B_n$
and $Y_n=(I_n+I_{n-1})\setminus (A_n+A_{n-1})$.
Let $\epsilon_n=1/n^2$.
Note that $a_n^2(y_{n-1}-x_{n-1})\geq 16\epsilon_n^{-2}\log(2y_n\epsilon_n^{-1})$ for every $n$ large enough.
By Corollary \ref{lm:proba}, for every $n$ large enough,
with probability at least $1-\epsilon_n$,
we have $|X_n\cap (2x_n+[0,y-1])|\leq \epsilon_n y$
and $|Y_n\cap (x_{n-1}+x_n+[0,y-1])|\leq \epsilon_n y$ for every $y\geq \epsilon x_n$.
Since $\sum_n \epsilon_n$ converges, by Borel--Cantelli, with probability 1,
we have $|X_n\cap (2x_n+[0,y-1])|\leq \epsilon_n y$
and $|Y_n\cap (T_{n-1}+T_n+[0,y-1])|\leq \epsilon_n y$ for every $y\geq \epsilon x_n$ and every $n$ large enough simultaneously.
Thus we can apply Corollary \ref{cor:dens0}:
we infer that the upper asymptotic density of $2B$ is almost surely just that of $2I$, which is
$b$ as was shown in the proof of Lemma \ref{lm:afa}.
In particular, there exists $B\subset\N$ such that $(\udens{B},\ldens{2B})=(a,b)$ as desired.
\end{proof}

\section{A diophantine intermezzo}
\label{sec:dioph}
In order to present our construction of a family of sets $A$ whose density profiles
$p_A$ cover a full open set, we need to state and prove  some diophantine results, which build upon those from \cite{BH}.
For a subset $A\subset\R/\Z$ of the circle, a set 
$I\subset\N$ of integers and $\lambda\in \R/\Z$, let 
$B_{\lambda,I,A}=\{n\in I :\lambda n\in A\}$.
When $I=\N$ we may omit it. Further, in a context where $\lambda$ is fixed, it may also be omitted.
Let $\mu$ denote the Lebesgue (unit) measure on $\T$.
\begin{theorem}[Weyl's equidistribution theorem]
\label{th:weyl}
    For any interval $A\subset\R/\Z$ of the unit circle and $\lambda\in \R\setminus\Q/\Z$,
    for any interval of integers
    $I\subset\N$,
    we have
    $$
   \abs*{ \{n\in I:n\lambda\in A\}}=(\mu(A)+o(1))|I|
    $$as $|I|$ tends to infinity.
    In particular, $B_{\lambda,A}$ has natural density $\alpha$.
\end{theorem}
It will turn out handy to make sure the error term in the equidistribution theorem decays reasonably rapidly.
For that, we recall a quantitative version of the equidistribution theorem due to Erd\H{o}s-Turan.
\begin{theorem}
 For any sequence $s_j$ of elements of the circle $\T=\R/\Z$ and any interval $A\subset\T$, we have for any integer interval $I$ of cardinality $n\in\N$ and integer $m$ the bound
\begin{equation}
    \label{eq:erdosturan}
\abs*{\frac{1}{\abs{I}}\abs*{\{ j \in I : s_j \in A\}} - \alpha}\ll
\frac{1}{m} + \frac{1}{n}\sum_{k=1}^m\frac{1}{k}\sum_{j=1}^n
e^{2i\pi s_j k}\end{equation}
where the implied constant is absolute.
\end{theorem}
Applying this with $s_j=\theta j$
for some irrational number $\theta$
and using the standard exponential sum bound
$\abs*{\sum_{j=1}^k
e^{2i\pi j\theta}}
\leq \frac{1}{2||\theta||}$
where $||\theta||=\min_{k\in\Z}|\theta-k|$, we obtain
\begin{equation}
    \label{eq:bound}
    \abs*{B_{\theta,I,A}-\alpha}\ll
    \frac{1}{m} + \frac{1}{n}\sum_{k=1}^m\frac{1}{k\nor{k\theta}}.
\end{equation}
Now the bound in equation \eqref{eq:bound} can be quite
poor for certain irrational numbers, but satisfactory for other. In particular, for algebraic numbers, by Roth's theorem,
we have $\nor{k\theta}\ll_\epsilon k^{-1-\epsilon}$
for any $\epsilon$.
Denote $\mu(A)=\alpha$.
In that case,
equation \eqref{eq:bound} becomes 
\begin{equation}
        \label{eq:boundAlg}
    \abs*{B_{\theta,I,A}-\alpha}\ll_\epsilon
    \frac{1}{m} + \frac{m^{1+\epsilon}}{n}.
\end{equation}
Choosing $m=n^{1/2}$,
we infer
\begin{equation}
        \label{eq:boundAlgEps}
    \abs*{B_{\theta,I,A}-\alpha}\ll_\epsilon
    n^{-1/2+\epsilon}.
\end{equation}
Selecting $\epsilon=1/6$,
one obtains
\begin{equation}
        \label{eq:boundAlgEps2}
    \abs*{B_{\theta,I,A}-\alpha}\ll
    n^{-1/3}.
\end{equation}
where the implied constant is absolute.
\begin{lemma}
\label{lm:diophSum}
Let $\lambda\in \R/\Z$ be irrational algebraic and  $A_1,A_2$ two open intervals of the circle $\T=\R/\Z$.
Let $I_1,I_2$ be two integer intervals.
Let $B_i=\{n\in I_i :\lambda n \in A_i\}$ for $i\in\{1,2\}$.
Let $S=\{n\in I_1+I_2 : \lambda n \in A_1+A_2\}$.
Then $B_1+B_2\subset S$. Conversely, 
for any $\eta\in (0,1)$ and $y\in\N$,
if $\min_{i\in\{1,2\}}|I_i|\geq C\eta^{-3}$
and $y\geq C\eta^{-4}$ for some absolute constant $C>0$, then
$R_{B_1,B_2}(x)\geq \eta^{-2}$ for all but $\eta y$ elements $x\in S\cap [\min(I_1+I_2)+y]$.
\end{lemma}
\begin{proof}
The inclusion $B_1+B_2\subset S$ holds since for any $n\in B_1,m\in B_2$,
we have $n+m\in I_1+I_2$ and $\lambda (n+m)=\lambda n+\lambda m\in A_1+A_2$.

We now prove the second claim of the Lemma.
For a set $X\subset\T$ and a positive real $\eta$, 
let $X^\eta$ be the set of $x\in X$ such that $(x-\eta,x+\eta)\subset X$.

For any $k\in \N$, 
denote by $k=a_i(k)+b_i(k)$ for $0\leq i<R_{I_1,I_2}(k)$ the restricted representations of $k$ with respect to $I_1,I_2$;  we may order them in such a way that $a_i(k)=a+i\in I_1$ and $b_i(k)=b-i\in I_2$, for each $i$ and some $a=a(k)$ and $b=b(k)=k-a$.
From Lemma \ref{lm:repres}, we know that $R(n)\geq f$ for all but $O(f)$ elements
of $I_1+I_2$.
Now
suppose $k\lambda \in (A_1+A_2)^{\eta/3}$.
Suppose $k=\ell+m$.
Observe that the conditions $\ell\lambda \in A_1$
and $m\lambda\in A_2$ are equivalent to $\ell\lambda\in A_1\cap (k\lambda-A_2)$.
Since $k\lambda\in (A_1+A_2)^{\eta/3}$, the intersection 
$A_1\cap (k\lambda-A_2)$ is a circle interval of length at least $\eta/3$.
Thus, the number of $\ell$ in an interval $I$ such that
$\ell\lambda\in A_1\cap (k\lambda-A_2)$ is at least $\abs{I}(\eta/3 -o(1))$ as $\abs{I}$ tends to infinity, by Weyl's equidistribution theorem (Theorem \ref{th:weyl}).
In particular, this number is at least $\eta\abs{I}/4$
as soon as $\abs{I}\geq C\eta^{-3}$ for some constant $C$  by equation \eqref{eq:boundAlgEps2} since $\lambda$ is irrational algebraic.
We may assume $C\geq 4$.
Therefore, when $k$ is such that
$R_{I_1,I_2}(k)\geq C\eta^{-3}$ and 
$\{k\lambda\}\in (A_1+A_2)^{\eta/3}$, we have
$R_{B_1,B_2}(k)\geq \eta R_{I_1,I_2}(k)/4$.
Now 
if $\min(|I_1|,|I_2|)\geq 2C\eta^{-3}$,
all but $X=O(\eta^{-3})$ elements $k\in I_1+I_2$
satisfy $R_{I_1,I_2}(k)\geq  C\eta^{-3}$.
If additionally $R_{B_1,B_2}(k)\geq \eta R_{I_1,I_2}(k)/3$, this yields $R_{B_1,B_2}(k)\geq C\eta^{-2}.$
Further,
$(A_1+A_2)\setminus (A_1+A_2)^{\eta/3}$ is a union of two intervals of length $\eta/3$,
so
if $y\geq C'\eta^{-3}$, for some constant $C'$
only
$3\eta y/4$ elements $k\in S\cap [\min(I_1+I_2)+y]$ satisfy $k\lambda\notin (A_1+A_2)^{\eta/3}$. 
If 
$y\geq C''\eta^{-4}$ for some absolute constant $C''$,
the term $X=O(\eta^{-3})$ is in turn less than $\eta y/4$ and we conclude, upon adjusting $C$ as required.
    \end{proof}

\section{A full dimensional polyhedron in $\cD$}
\label{sec:quad}
In this section, we 
prove Theorem \ref{th:quadruples}, whose statement we recall here.
\quadruples*

The subset of $(0,1]^4$ given above contains the non-empty interior of a polyhedron. In particular it has a positive Lebesgue measure.
We present a construction based on a combination of all the diophantine, probabilistic, and archimedean methods.
\begin{proof}
    Let $\alpha_1,\alpha_2\in (0,1]^2$ and assume $\alpha_1\leq \alpha_2$.
    Let $\beta_1,\beta_2$ be two real numbers in $(0,1]$ and assume $\beta_1 \alpha_1\leq \beta_2 \alpha_2$.
    Our aim is to show that there exists a set whose density profile is
    $(\alpha_1\beta_1,\alpha_2\beta_2,\min(1,\alpha_1+\alpha_2),\min(1,2\alpha_2))$.
    
    Let $A_i$ be the interval $[0,\alpha_i]$ in the circle $\R/\Z$, for $i=1,2$.
    We extend these definitions to $\alpha_k,\beta_k,A_k$ for all $k\in\N$ in a 2-periodic way.
    Consider an increasing sequence $(a_n)_{n\in\N}$ of positive integers satisfying $a_n/a_{n+1}\rightarrow 0$. It will be convenient to assume $a_n$ does not grow to fast either;
    for definiteness we set $a_n=2^{n^2}$.
    Then let $I_n=[a_{n-1},a_n]$. Let $\theta\in (0,1)$ be irrational algebraic.
    Let $B_k=\{n\in I_k :n\theta\in A_k\}$,
    where $\alpha_k$ is defined by reducing
    $k$ modulo 2, for every $k\in\N$.
   Let $C_k$ a Bernoulli subset of $B_k$ of parameter $\beta_k$,
    where $\beta_k$ is defined by reducing
    $k$ modulo 2, for every $k\in\N$.
    Finally, let $C=\bigcup_n C_n$.

\textbf{Claim}   Almost surely, $\ldens{C}=\beta_1 \alpha_1$
    and  $\udens{C}=\beta_2 \alpha_2$.  

\textbf{Proof of the Claim}
First, we have 
$\abs{C\cap [a_n]}=O(a_{n-1})+\abs{C_{n}}$.
Now $\abs{C_{n}}=\sum_{k\in B_{n}}1_{k\in C_{n}}$.
By the law of large numbers, almost surely as $n$ tends to infinity, $\sum_{k\in B_{n}}1_{k\in C_{n}}=(\beta_{n}+o(1))\abs{B_{n}}$.
By Weyl's equidistribution theorem (Theorem \ref{th:weyl}), $\abs{B_{n}}=\abs{I_{n}}(\alpha_n+o(1))$ as $n$ tends
to infinity.
Using the fact that $a_{n-1}/a_n$ tends to 0,
we infer $\ldens{A}\leq\liminf_n \abs{C\cap [a_n]}/a_n=\alpha_1\beta_1$
and $\udens{A}\geq \limsup_n \abs{C\cap [a_n]}/a_n=\alpha_2\beta_2$.

To conclude, it suffices to show that
for every $\epsilon >0$, with probability 1, for any $x\in\N$ large enough (in terms of $\epsilon$),
we have $|C\cap [x]|/x\in [\alpha_1\beta_1-\epsilon,\alpha_2\beta_2+\epsilon]$.
Let $\epsilon>0$. Let $x\in\N$ and $n$ such that $x\in I_{n+1}$.
Write $x=a_{n}+y$ for some $y\geq 0$.
     We have $|C\cap [x]|=|C\cap [a_{n}]| +|C\cap (a_{n},a_{n}+y]|$,
    hence
    \begin{equation}
        \label{eq:barycentren}
        \frac{|C\cap [x]|}{x}=\frac{a_{n}}{x}\frac{|C\cap [a_{n}]|}{a_{n}}
+\frac{y}{x}\frac{|C\cap (a_{n},a_{n}+y]|}{y}
    \end{equation}
We have $\frac{\abs{|C\cap [a_{n}]|-\alpha_n\beta_n a_{n}}}{a_{n}}\leq \epsilon/2$, by the considerations above, when $n$ is large enough -- hence also when $x$ is large enough -- with
probability 1.
Either $y\leq \epsilon x/2$, in which case 
$$
\frac{y}{x}\frac{|C\cap (a_{n},a_{n}+y]|}{y}\in [0,\epsilon/2]
$$ whence
we conclude invoking equation \eqref{eq:barycentren}
that  $|C\cap [x]|/x\in [\alpha_1\beta_1-\epsilon/,\alpha_2\beta_2+\epsilon]$, for all $x$ large enough (of this type) with probability 1.
Now suppose $y>\epsilon x/2$.
Let $w=|C\cap (a_{n},a_{n}+y]|$ and $z=|B\cap (a_{n},a_{n}+y]|$. These quantities depend
on $y$ and therefore on $x$.
Write $\beta=\beta_{n+1},\alpha=\alpha_{n+1}$ for simplicity.
Then $\Pr(\abs{w-\beta z}\geq \eta z)\leq \exp(-2z\eta^2)$ by Hoeffding's bound, for any $\eta >0$.
We will take $\eta=\epsilon/4$
Also $|z-\alpha y|\leq \epsilon y/4$, by Weyl's theorem if we assume $n$ large enough, since $y>\epsilon a_n/2$.
By the union bound
the probability that there exists $x\in (a_{n},a_{n+1}]|$ such that
$\abs{w-\beta z}\geq \eta z$ is at most
$a_{n+1}\exp(-ca_n)$ where $c=2(\alpha-\epsilon/2)\eta^2>0$.
Now since $a_{n+1}\leq \exp(ca_n/2)$ for $n$ large enough,
we infer that the series $\sum_n a_{n+1}\exp(-ca_n)$ converges; by
the Borel-Cantelli lemma,
with probability 1, for all $n$ large enough and all $w\in ( (1+\epsilon/2)a_n,a_{n+1}]$, we have
$$\abs{\abs{C\cap (a_{n},x]}-\beta\alpha y}\leq \abs{\abs{C\cap (a_{n},x]}-\alpha\abs{B\cap (a_{n},x]}}+\alpha \abs{\abs{B\cap (a_{n},x]}-\beta y}
\leq \epsilon/2.
$$
Injecting this result in \eqref{eq:barycentren},
we conclude the proof of the Claim.

    Now we determine the upper and lower asymptotic densities of $2C$.
    We show that they are $\alpha_1+\alpha_2$ and $2\alpha_2$, respectively.
We decompose the proof into several claims.

\textbf{Claim}
$\udens{2C}\leq \min(1,2\alpha_2)$

This claim follows from the fact that $2C\subset B_{\theta,\N,2A_2}=\{n\in\N : n\theta\in 2A_2\}$ and Weyl's equidistribution theorem since $2A_2$ has Lebesgue measure $\min(1,2\alpha_2)$.

\textbf{Claim}
Almost surely, 
$B_{\theta,\N,A_1+A_2}\setminus 2C$
has density 0. In particular $\ldens{2C}\geq \min(1,\alpha_1+\alpha_2)$.

To prove this claim, we first observe that
$2C\supset \bigcup_n C_n+C_{n+1}$.
Let $S_{n}=B_{I_n+I_{n+1},A_1+A_2}$.
Let $X_n=S_n\setminus C_n+C_{n+1}$ and
$X=\bigcup_n X_n=B_{\theta,\N,A_1+A_2}\setminus 2C$. Our goal is to prove that almost surely,
$X$ has density 0.
Note that $X_n\subset I_n+I_{n+1}=\llbracket b_n,b_{n+1}\rrbracket $
where $b_n=a_n+a_{n-1}$.
Let $\epsilon>0$. We will prove that $\udens{X}\leq \epsilon$ almost surely, which will conclude.
Let $x\in\N$. We shall try to bound $|X\cap [x]|$.
Let $n$ be such that $b_n\lt x \leq b_{n+1}$.
We have
$$
|X\cap [x]|=\sum_{k\lt n}\abs{X_k}+\abs{X_n\cap (b_n,x\rrbracket}.
$$
We write $x=b_n+y$.
As in equation \eqref{eq:barycentren}, we may write
\begin{equation}
        \label{eq:barycentre2}
        \frac{|X\cap [x]|}{x}=\frac{b_{n}}{x}\frac{|X\cap [b_{n}]|}{b_{n}}
+\frac{y}{x}\frac{|X\cap (b_{n},b_{n}+y]|}{y}
    \end{equation}
If $y<\epsilon b_n/2$, then
$\frac{|X\cap [x]|}{x}\leq \frac{|X\cap [b_{n}]|}{b_{n}}+\epsilon/2$.
So
we may assume that $x=b_n+y$ with $y\geq \epsilon b_n$.
Let $\eta \mt 0$ to be determined later.
When $n$ is large enough, $\epsilon b_n$ is larger than $C\eta^{-4}$ so
we can apply Lemma \ref{lm:diophSum};
we infer that there exists a positive constant $C$ such that 
    there exists a set $K\subset S_n\cap (b_{n},b_{n}+y]$ such that 
     $\abs{S_n\cap (b_{n},b_{n}+y]\setminus K}\leq\eta y$
    and each element of $K$ admits at least $2\eta^{-2}$ representations
    $b_1+b_2$ with $(b_1,b_2)\in B_n\times B_{n+1}$.
    By Lemma \ref{lm:basicProba}, this conclusion implies that
    $\abs{K\setminus (C_n+C_{n-1})}\leq \abs{K}e^{-\eta^{-2}\beta_1\beta_2}$ 
    with probability at least $1-e^{-\eta^{-2}\beta_1\beta_2}$.
    Therefore, assuming $\abs{I_n}\geq C\eta^{-4}$, we have 
    $\abs{X\cap (b_{n},b_{n}+y]}\leq y(\eta+e^{-\eta^{-2}\beta_1\beta_2})$
    with probability
    at least $1-e^{-\eta^{-2}\beta_1\beta_2}$.
    By the union bound,
    with probability at least
    $1-b_{n+1}e^{-\eta^{-2}\beta_1\beta_2}$,
    for all $x\in ((1+\epsilon)b_n,b_{n+1}\rrbracket$
    we have $\abs{X\cap (b_{n},x]}\leq y(\eta+e^{-\eta^{-2}\beta_1\beta_2})$.
    Now we select $\eta=\eta_n=1/n^2$; thus $(\epsilon b_n\geq \epsilon 2^{(n-1)^2}\geq C\eta^{-4}$ for
    $n$ large enough, and $\sum_n b_{n+1}e^{-\eta_n^{-2}\beta_1\beta_2}<\infty$.
    By Borel-Cantelli, with probability 1, for all $n$ large enough and all $x\in ((1+\epsilon)b_n,b_{n+1}\rrbracket$,
    $\abs{X\cap (b_n,x]}\leq \epsilon (x-b_n)$.
    We conclude that $\Pr(\dens{X}\leq\epsilon)=1$,
    and this being true for any $\epsilon\mt 0$,
    we infer that 
    $\Pr(\dens{X}=0)=1$.
    That is, $\Pr(\dens{B_{\theta,\N,A_1+A_2}}\setminus 2C=0)=1$ and hence $\ldens{2C}\geq \dens{B_{\theta,\N,A_1+A_2}}=\min(1,\alpha_1+\alpha_2)$ almost surely, by the equidistribution theorem.

    \textbf{Claim}
    For all $\epsilon\mt 0$, 
    $\abs{2C\cap [a_n]}\leq a_n(\mu(A_1+A_2)+\epsilon)$
    for all odd $n$ large enough.
In particular, $\ldens{2C}\leq \min(1,\alpha_1+\alpha_2)$.
    \begin{proof}[Proof of Claim]
        We have $2C\cap [a_n]\subset (B_{n}+\bigcup_{k\leq n}B_k)\cap [a_n]\cup [a_{n-1}]$.
        Now $(B_{n-1}+\bigcup_{k\lt n}B_k)\cap [a_n]\subset B_{[a_n],A_1+A_2}$ since $n$ is odd.
        Now  $\abs{B_{[a_n],A_1+A_2}}\leq (\alpha_1+\alpha_2+\epsilon/2)a_n$ when $n$ is large enough by the equidistribution theorem \ref{th:weyl}.
        Of course, $a_{n-1}\leq a_n\epsilon/2$ for $n$ large enough too.
        This proves the claim.
    \end{proof}

    \textbf{Claim}
    For every $\epsilon \mt 0$,
    when $n$ is even and large enough, almost surely,
    $\abs{2C\cap [2a_n]}\geq 2a_n(\mu(2A_2)-\epsilon)$.
    In particular, almost surely,
    $\udens{2C}\geq \min(1,2\alpha_2)$

    \begin{proof}
        We have
        $2C\cap [2a_n]\supset C_n+C_n$.
        Take $\eta=n^{-1}$, so that $\abs{C_n}\geq C\eta^{-4}$
        for any constant $C$ for $n$ large enough.
        By Lemma \ref{lm:diophSum},
        all but $\eta a_n$ elements
        of $B_{2I_n,2A_2}$ admit at least $2\eta^{-2}$
        representations $b_1+b_2$ with $b_1\lt b_2\in B_n\times B_n$.
        By Lemma \ref{lm:basicProba},
        it follows that $\abs{C_n+C_n}\geq \abs{B_{2I_n,2A_n}}(1-\eta-e^{-\beta_2^2\eta^{-2}})$
        with probability $1-e^{-\beta_2^2\eta^{-2}}$.
        Applying Borel-Cantelli,
        we infer that almost surely
        $\abs{C_n+C_n}\geq 2a_n(2\alpha_2-\epsilon)$
        for $n$ large enough. This concludes.
    \end{proof}
    
    The combination of the five claims yields that almost surely
    $p_C=G(\alpha_1,\alpha_2,\beta_1,\beta_2)$
    where
    $G(\alpha_1,\alpha_2,\beta_1,\beta_2)=(\beta_1\alpha_1,\beta_2\alpha_2,\min(1,\alpha_1+\alpha_2),\min(1,2\alpha_2))$,
    in particular there exists a set which has this density profile.
    Observe that for any  $(\alpha_1,\alpha_2,\beta_1,\beta_2)\in (0,1]^4$ 
  satisfying
  $\alpha_1\leq\alpha_2$ and  $\beta_1\alpha_1\leq\beta_2\alpha_2$, the quadruplet
  $G(\alpha_1,\alpha_2,\beta_1,\beta_2)\in (0,1]^4$ satisfies the system of
  equations given in Theorem \ref{th:quadruples}.
  
  To conclude the proof, let $a,b,c,d$ be as in the statement of Theorem \ref{th:quadruples}.
  Then we find a quadruplet $\sigma:=(\alpha_1,\alpha_2,\beta_1,\beta_2)\in (0,1]^4$ 
  satisfying
  $\alpha_1\leq\alpha_2$ and  $\beta_1\alpha_1\leq\beta_2\alpha_2$
  such that
  $(a,b,c,d)=G(\sigma)$.
  It is easy to control that $\sigma$ defined below is such a quadruplet.
  \begin{itemize}
      \item If $d=c=1$,
      take $\sigma=(a,b,1,1)$.
      \item If $d=1>c$, take $\sigma=(c-b,b,a/(c-b),1)$.
      \item If $d<1$, take $\sigma=(c-d/2,d/2,a/(c-d/2),2b/d)$
  \end{itemize}
\end{proof}

    \section{Application to certain three-dimensional projections of $\cD$
    }
    \label{sec:3dim}
Hegyv\'ari, Hennecart and Pach 
\cite[Proposition 2.2]{hennecartAustralasian}
proved that there is a set $A\subset\N$ which has a density
whereas $2A$ does not. They asked for the possible values of the triplet $(\dens{A},\ldens{2A},\udens{2A})$, providing only
$(3/7,6/7,1)$ by an explicit construction.
This amounts to asking for the projection of $\cD$ onto the three dimensional space $\{(a,b,c,d)\in \R^4:a=b\}$.
We supply here a partial solution to the question, comprising as a particular case
the triplet $(3/7,6/7,1)$ already found (by a different method).
\begin{cor}
    The triplet $(\dens{A},\ldens{2A},\udens{2A})$
    may take any value $(a,c,d)\in [0,1]^4$
    satisfying the system
    $$
    \left\{
\begin{array}{ccc}
    c & \leq & d \\
   \min(2a,1)&\leq   & c \\
   d&< & 2c\\
\end{array}
    \right.
    $$
\end{cor}
This corollary is an obvious consequence of Theorem \ref{th:quadruples}.

A set $A$ which has no natural density such that $2A$ does is  more natural to think of and the aforementioned authors also asked for the possible values of the triplet
$(\ldens{A},\udens{A},\dens{2A})$.
Again, as a direct application of Theorem \ref{th:quadruples}, we are able to give a partial answer to their question.
\begin{cor}
    The triplet $(\ldens{A},\udens{A},\dens{2A})$
    may take any value $(a,b,c)\in [0,1]^3$
    satisfying the system of equations 
     $$\left\{
    \begin{array}{ccc}
        a    & \leq  & b \\
        2b &\leq  &c \\
    \end{array}
    \right.
   $$
\end{cor}

\section{Further remarks and opening}
\label{sec:remarks}
Theorem \ref{th:quadruples} is certainly far from optimal.
First, we could relax the strict inequalities $a>0,b>0,2c>d$, at the cost of replacing, in the diophantine construction, the parameters $\alpha_1,\alpha_2$ by sequences tending to infinity slowly enough. We found this would complicate seriously the proof while improving only marginally the result and so did not include it.

The main clearly suboptimal feature of our construction is that it does not allow $\ldens{2A}$ to be small and $\udens{A}$ to be large simultaneously, whereas we know from Theorem \ref{th:23}
that this can happen.

However, it is clear that there are constraints on the quadruplet $p_A$ 
that are not reflected in its  subpairs.
For instance, we know that if $A$ is reduced (i.e. $\gcd(A-A)=1$)
and $\udens{2A}$ is as small as possible given $\udens{A}$, i.e. $\udens{2A}=f(\udens{A})$,
where $f$ is as above the function $[0,1]\rightarrow [0,1]$ defined by $f(\alpha)=\min(3\alpha/2,(1+\alpha)/2)$,
then $\ldens{A}=\ldens{2A}=0$, unless there is an non-archimedean reason, i.e. $2A$ is eventually periodic. This follows from Jing's inverse theorem \cite{jin,jinSol}.
Bordes \cite{bordes} showed that such constraints occur even a little beyond the regime $\udens{2A}=f(\udens{A})$.
Thus there is a constant $\alpha_0\in (0,1)$ such that whenever $A$ is reduced, $\udens{A}<\alpha_0$
and $\udens{2A}<(3/2+\epsilon)\udens{A}$, then $\ldens{A}<2\epsilon(1-\udens{A})^{-1}+O(\epsilon^2)$, for any $\epsilon$ small enough.
In view of such results, which are not sharp either, it seems very difficult to guess,
let alone determine, the exact value of the set $\cD$.

\end{document}